\documentclass[a4paper,11pt]{amsart}
\usepackage{graphicx}
\usepackage{mathptmx}
\usepackage{amsmath}
\usepackage{amssymb}
\usepackage{enumitem}
\usepackage{xcolor}

\newmuskip\pFqmuskip

\newcommand*\pFq[6][8]{%
  \begingroup % only local assignments
  \pFqmuskip=#1mu\relax
  \mathcode`=\string"8000
  \begingroup\lccode`\~=`\,
  \lowercase{\endgroup\let~}\pFqcomma
  F^{#2}_{#3}{\left(\genfrac..{0pt}{}{#4}{#5}\bigg|#6\right)}%
  \endgroup
}
\newcommand{\pFqcomma}{\mskip\pFqmuskip}

\newtheorem{theorem}{Theorem}

\newtheorem{remark}[theorem]{Remark}

\begin{document}

\title[]{Normal ordering associated with $\lambda$-Stirling numbers\\ in $\lambda$-Shift algebra}

\author{Taekyun Kim $^{1,*}$}
\address{Department of Mathematics, Kwangwoon University, Seoul 139-701, Republic of Korea}
\email{tkkim@kw.ac.kr}

\author{Dae San Kim $^{2}$}
\address{Department of Mathematics, Sogang University, Seoul 121-742, Republic of Korea}
\email{dskim@sogang.ac.kr}

%\author{Hye Kyung Kim $^{3}$}
%\address{Department of Mathematics Education, Daegu Catholic University, Gyeongsan 38430, Republic of Korea}
%\email{hkkim@cu.ac.kr}

\subjclass[2020]{11B73; 11B83}
\keywords{$\lambda$-shift algebra; normal ordering; unsigned $\lambda$- Stirling numbers of the first kind; $\lambda$-$r$-Stirling numbers of the first kind}
\thanks{* is corresponding author}

\begin{abstract}
The Stirling numbers of the second kind are related to normal orderings in the Weyl algebra, while the unsigned Stirling numbers of the first kind are related to normal orderings in the shift algebra. Kim-Kim introduced a $\lambda$-analogue of the unsigned Stirling numbers of the first kind and that of the $r$-Stirling numbers of the first kind. In this paper, we introduce a $\lambda$-analogue of the shift algebra (called $\lambda$-shift algebra) and investigate normal orderings in the $\lambda$-shift algebra. From the normal orderings in the $\lambda$-shift algebra, we derive some identities about the $\lambda$-analogue of the unsigned Stirling numbers of the first kind .

\end{abstract}

 \maketitle

\markboth{\centerline{\scriptsize  Normal ordering associated with $\lambda$-Stirling numbers in $\lambda$-Shift algebra}}
{\centerline{\scriptsize   T. Kim, D. S. Kim}}

\section{Introduction}
%2
The Stirling number of the first kind $S_{1}(n,k)$ is defined in such a way that the unsigned Stirling number of the first kind ${n \brack k}=(-1)^{n-k}S_{1}(n,k)$ enumerates the number of permutations of the set $[n]=\{1,2,3,\dots,n\}$ which are products of $r$ disjoint cycles. \\
The unsigned $r$-Stirling number of the first kind ${n \brack k}_{r}$ is the number of permutations of $[n]$ with exactly $k$ disjoint cycles in such a way that the numbers $1,2,\dots,r$ are in distinct cycles. \par
In \cite{A}, introduced are a $\lambda$-analogue of the unsigned Stirling numbers of the first kind ${n \brack k}_{\lambda}$ and that of the unsigned $r$-Stirling numbers of the first kind ${n \brack k}_{r,\lambda}$ respectively as a $\lambda$-analogue of ${n \brack k}$ and that of ${n \brack k}_{r}$,\,\, (see \eqref{eq08}, \eqref{eq09}). \par
The Stirling numbers of the second kind appear as the coefficients in normal orderings in the Weyl algebra (see \eqref{eq10}, \eqref{eq11}), while the unsigned Stirling numbers of the first kind appear as those in normal orderings in the shift algebra $S$ (see \eqref{eq12}, \eqref{eq13}). \par
The aim of this paper is to introduce the $\lambda$-shift algebra $S_{\lambda}$ (for any $\lambda \in \mathbb{C}$), which is a $\lambda$-analogue of $S$ (see \eqref{eq14}), and to investigate normal orderings in the $\lambda$-shift algebra. In addition, from the normal orderings in the $\lambda$-shift algebra $S_{\lambda}$, we derive some identities about the unsigned $\lambda$-Stirling numbers of the first kind. \par

The outline of this paper is as follows. In Section 1, we recall the $\lambda$-falling factorial numbers, the falling factorial numbers, the $\lambda$-rising factorial numbers and rising factorial numbers. We remind the reader of the unsigned $\lambda$-Stirling numbers of the first kind and the $\lambda$-$r$-Stirling numbers of the first kind. We recall the Weyl algebra and the normal ordering result in that algebra. We remind the reader of the shift algebra and the normal ordering result in that algebra. Finally, we define the $\lambda$-shift algebra as a $\lambda$-analogue of the shift algebra. Section 2 is the main result of this paper. We derive normal ordering results in $S_{\lambda}$ in Theorem 1 and Theorem 2 where ${n \brack k}_{\lambda}$ and ${n+r \brack k+r}_{r, \lambda}$ appear respectively as their coefficients. We obtain three other normal ordering results in Theorem 3. In Theorem 4, we get a recurrence relation for the unsigned $\lambda$-Stirling numbers of the first kind. In Theorem 6, we get another expression of the defining equation in \eqref{eq08} in terms of the $\lambda$-shift operator (see \eqref{eq31}). In Theorem 7, we show a $\lambda$- analogue of the dual to Spivey's identity (see Remark 8). Finally, we conclude this paper in Section 3. For the rest of section, we recall what are needed throughout this paper. \par

\vspace{0.1in}

For any $\lambda \in \mathbb{C}$, the $\lambda$-falling factorial sequence is defined by
\begin{equation}\label{eq01}
\begin{split}
(x)_{0,\lambda}=1, \ \  (x)_{n,\lambda}=x(x-\lambda)\cdots(x-(n-1)\lambda), \quad (n\geq 1), \quad (\text{see \cite{4,11}}).
\end{split}
\end{equation}
In particular, the falling factorial sequence is given by
\begin{equation}\label{eq02}
\begin{split}
(x)_{0}=1, \ \  (x)_{n}=x(x-1)\cdots(x-(n-1)), \quad (n\geq 1).
\end{split}
\end{equation}
Note that $\lim_{\lambda\rightarrow 1}(x)_{n,\lambda}=(x)_{n}$.

For any $\lambda \in \mathbb{C}$, the $\lambda$-rising factorial sequence is defined by
\begin{equation}\label{eq03}
\begin{split}
\langle x\rangle_{0,\lambda}=1, \,\, \langle x\rangle_{n,\lambda}=x(x+\lambda)\cdots(x+(n-1)\lambda), \quad (n\geq 1), \quad (\text{see \cite{4,11}}).
\end{split}
\end{equation}
Especially, the rising factorial sequence is given by
\begin{equation}\label{eq04}
\begin{split}
\langle x\rangle_{0}=1, \,\,  \langle x \rangle_{n}=x(x+1)\cdots(x+(n-1)), \quad (n\geq 1).
\end{split}
\end{equation}
Observe that $\lim_{\lambda\rightarrow 1} \langle x\rangle_{n,\lambda}=\langle x\rangle_{n}$.

With the notation in \eqref{eq02}, the Stirling numbers of the first kind are defined by
\begin{equation}\label{eq05}
\begin{split}
(x)_n=\sum_{k=0}^n S_1(n,k)x^k,\quad (n\geq 0), \quad (\text{see \cite{5,6,19,23}}).
\end{split}
\end{equation}
In addition, the unsigned Stirling numbers of the first kind are given by ${n \brack k}=(-1)^{n-k}S_1(n,k), \,\, (n,\ k\geq0)$. \\
The Stirling numbers of the second kind are defined by
\begin{equation}\label{eq06}
\begin{split}
x^n=\sum_{k=0}^n {n \brace k}(x)_k, \quad (n \geq 0).
\end{split}
\end{equation} \par
Recently, with the notation in \eqref{eq01} the $\lambda$-Stirling numbers of the first kind, which are $\lambda$-analogues of the Stirling numbers of the first kind, are defined by
\begin{equation}\label{eq07}
\begin{split}
(x)_{n,\lambda}=\sum_{k=0}^n S_{1,\lambda}(n,k)x^k, \quad (n\geq0), \quad  (\text{see \cite{11}}).
\end{split}
\end{equation}
In addition, with the notation in \eqref{eq03} the unsigned $\lambda$-Stirling numbers of the first kind are defined by
%3
\begin{equation}\label{eq08}
\begin{split}
\langle x\rangle_{n,\lambda}=\sum_{k=0}^n {n \brack k}_\lambda x^k, \quad (n\geq0), \quad(\text{see \cite{11}}).
\end{split}
\end{equation}
Note that $\lim_{\lambda\rightarrow 1}S_{1,\lambda}(n,k)=S_1(n,k)\,\, (\mathrm{see} \eqref{eq05}), \,\, \lim_{\lambda\rightarrow 1}{n \brack k}_\lambda={n \brack k}$. \par

For $r \in \mathbb{N}\cup\{0\}$, the $\lambda$-$r$-Stirling numbers of the first kind, which are $\lambda$-analogues of the $r$-Stirling numbers of the first kind, are defined by
\begin{equation}\label{eq09}
\begin{split}
\langle x+r\rangle_{n,\lambda}=\sum_{k=0}^n {n+r \brack k+r}_{r,\lambda}x^k, \quad (n \geq0), \quad (\text{see \cite{11}}).
\end{split}
\end{equation}

Note that $\lim_{\lambda\rightarrow 1}{n+r \brack k+r}_{r,\lambda}={n+r \brack k+r}$, where ${n+r \brack k+r}$ are the $r$-Stirling numbers of the first kind which are introduced by Broder  \ (\text{see \cite{3}}) and given by (see \eqref{eq04})

\begin{equation*}
\begin{split}
\langle x+r\rangle_{n}=\sum_{k=0}^n {n+r \brack k+r}_{r}x^k, \quad (n \geq0).
\end{split}
\end{equation*}

The Weyl algebra is the unital algebra generated by letters $a$ and $a^\dag$ satisfying the commutation

\begin{equation}\label{eq10}
\begin{split}
aa^\dag - a^\dag a=1, \quad (\text{see [1, 2, 7-21]}).
\end{split}
\end{equation}

Katriel proved that the normal ordering in Weyl algebra is given by (see \eqref{eq06})
\begin{equation}\label{eq11}
\begin{split}
(a^\dag a)^n=\sum_{k=0}^n{n \brace k}(a^\dag)^ka^k, \quad (\text{see [7-9]}).
\end{split}
\end{equation}
%4
From the definition of the Stirling numbers of the second kind and \eqref{eq11}, we note that
\begin{equation*}
\begin{split}
(a^\dag)^na^n=(a^\dag a)_n=a^\dag a(a^\dag a-1)\cdots(a^\dag a-n+1), \quad (n\geq 1).
\end{split}
\end{equation*}

The shift algebra $S$ is defined as the complex unital algebra generated by $a$ and $a^\dag$ satisfying the commutation relation

\begin{equation}\label{eq12}
\begin{split}
aa^\dag - a^\dag a=a, \quad (\text{see \cite{20}}).
\end{split}
\end{equation}
A word in $S$ is said to be in normal ordered form if all letters $a$ stand to the right of all letters $a^\dag$.

From \eqref{eq12}, we note that the normal ordering in the shift algebra $S$ is given by (see \eqref{eq11})
\begin{equation}\label{eq13}
\begin{split}
(a^\dag a)^n=\sum_{k=0}^n {n \brack k}(a^\dag)^k a^n, \quad (\text{see \cite{20}}).
\end{split}
\end{equation}
%5
For any $\lambda \in \mathbb{C}$, we consider a $\lambda$-analogue of the shift algebra $S$ which is defined as the complex unital algebra generated by $a$ and $a^\dag$ satisfying the commutation relation (see \eqref{eq12})
\begin{equation}\label{eq14}
\begin{split}
aa^\dag -a^\dag a=\lambda a.
\end{split}
\end{equation}
The $\lambda$-analogue of the shift algebra $S$ is called the $\lambda$-shift algebra and denoted by $S_{\lambda}$.

\medskip

\section{$\lambda$-analogues of normal ordering in the $\lambda$-shift algebra.}

%6
Let $S_{\lambda}$ be the  $\lambda$-shift algebra defined in \eqref{eq14}.
A word in $S_\lambda$ is said to be in normal ordered form if all letters $a$ stand to the right of all letters $a^\dag$.

In $S_\lambda$, by \eqref{eq14}, we get
\begin{equation*}
\begin{split}
(a^\dag a)^2&=(a^\dag a)(a^\dag a)=a^\dag(aa^\dag)a=a^\dag(\lambda a+a^\dag a)a \\
&=a^\dag(\lambda+a^\dag)a^2=\langle a^\dag \rangle_{2,\lambda}a^2,
\end{split}
\end{equation*}
\begin{equation*}
\begin{split}
(a^\dag a)^3&=(a^\dag a)(a^\dag a)(a^\dag a)=a^\dag(aa^\dag)(aa^\dag)a \\
&=a^\dag(\lambda+a^\dag)a(\lambda+a^\dag)a^2=a^\dag(\lambda+a^\dag)(a\lambda+aa^\dag)a^2 \\
&=a^\dag(\lambda+a^\dag)(2\lambda a+a^\dag a)a^2= a^\dag(a^\dag+\lambda)(a^\dag+2\lambda)a^3 \\
&=\langle a^\dag \rangle_{3,\lambda}a^3.
\end{split}
\end{equation*}

Continuing this process, we have
\begin{equation}\label{eq16}
\begin{split}
(a^\dag a)^n=\langle a^\dag \rangle_{n,\lambda}a^n, \quad (n\geq 1).
\end{split}
\end{equation}

%7
Thus, by \eqref{eq08} and \eqref{eq16}, we get (see \eqref{eq07}, \eqref{eq11}, \eqref{eq13})
\begin{equation}\label{eq17}
\begin{split}
(a^\dag a)^n =\sum_{k=0}^n {n \brack k}_\lambda (a^\dag)^k a^n.
\end{split}
\end{equation}

Therefore, by \eqref{eq17}, we obtain the following theorem.

\begin{theorem}
In $S_\lambda$, the unsigned $\lambda$-Stirling numbers of the first kind appear as the coefficients of $(a^\dag a)^n$ in normal ordered form, as it is given by
\begin{equation*}
\begin{split}
(a^\dag a)^n=\sum_{k=0}^n {n \brack k}_\lambda (a^\dag)^k a^n.
\end{split}
\end{equation*}
\end{theorem}

\medskip

For $r \geq 0$, by \eqref{eq14},  we get
\begin{equation*}
\begin{split}
((a^\dag+r)a)^2&=((a^\dag+r)a)((a^\dag+r)a)=(a^\dag+r)(aa^\dag+ra)a \\
&=(a^\dag+r)(a^\dag a+\lambda a+ra)a=(a^\dag+r)(a^\dag+r+\lambda)a^2 \\
&=\langle a^\dag+r\rangle_{2,\lambda}a^2,
\end{split}
\end{equation*}

and

\begin{equation*}
\begin{split}
((a^\dag+r)a)^3&=((a^\dag+r)a)((a^\dag+r)a)((a^\dag+r)a) \\
&=(a^\dag+r)a(a^\dag+r)(\lambda+r+a^\dag)a^2 \\
&=(a^\dag+r)(\lambda+r+a^\dag)a(\lambda+r+a^\dag)a^2 \\
&=(a^\dag+r)(a^\dag+r+\lambda)(\lambda a+ra+aa^\dag)a^2 \\
&=(a^\dag+r)(a^\dag+r+\lambda)(a^\dag+r+2\lambda)a^3=\langle a^\dag+r\rangle_{3,\lambda}a^3.
\end{split}
\end{equation*}

%8
Continuing this process, we have
\begin{equation}\label{eq18}
\begin{split}
((a^\dag+r))a^n=\langle a^\dag+r \rangle_{n,\lambda}a^n, \quad (n\geq 1).
\end{split}
\end{equation}

From \eqref{eq09} and \eqref{eq18}, we get
\begin{equation}\label{eq19}
\begin{split}
((a^\dag+r)a)^n&=\langle a^\dag+r \rangle_{n,\lambda}a^n \\
&=\sum_{k=0}^n{n+r \brack k+r}_{r,\lambda}(a^\dag)^ka^n.
\end{split}
\end{equation}

Therefore, by \eqref{eq19}, we obtain the following theorem.

\begin{theorem}

Let $r$ be a nonnegative integer. In $S_{\lambda}$, the $\lambda$-$r$-Stirling numbers of the first kind appear as the coefficients of $((a^\dag+r)a)^n$ in normal ordered form, as it is given by
\begin{equation*}
\begin{split}
((a^\dag+r)a)^n=\sum_{k=0}^n {n+r \brack k+r}_{r,\lambda}(a^\dag)^ka^n.
\end{split}
\end{equation*}

\end{theorem}

\medskip

%9
From \eqref{eq14}, we note that
\begin{equation}\label{eq20}
\begin{split}
a^ma^\dag&=a^{m-1}(aa^\dag)=a^{m-1}(a^\dag+\lambda)a\\
&=(a^{m-1}a^\dag)a+\lambda a^m=a^{m-2}(a^\dag a+\lambda a)a+\lambda a^m \\
&=(a^{m-2}a^\dag)a^2+2\lambda a^m=\cdots=(a^\dag+m\lambda)a^m,
\end{split}
\end{equation}

and

\begin{equation}\label{eq21}
\begin{split}
a(a^\dag)^n&=(aa^\dag)(a^\dag)^{n-1}=(\lambda+a^\dag)a(a^\dag)^{n-1} \\
&=(\lambda+a^\dag)(aa^\dag)(a^\dag)^{n-2}=(\lambda+a^\dag)^2a(a^\dag)^{n-2}=\cdots \\
&=(\lambda+a^\dag)^{n-1}aa^\dag=(\lambda+a^\dag)^n a.
\end{split}
\end{equation}

By \eqref{eq20}, we get
\begin{equation}\label{eq22}
\begin{split}
a^m(a^\dag)^n&=(a^ma^\dag)(a^\dag)^{n-1}=(a^\dag +m\lambda)a^m(a^\dag)^{n-1} \\
&=(a^\dag+m\lambda)(a^ma^\dag)(a^\dag)^{n-2}=(a^\dag+m\lambda)(a^\dag+m\lambda)a^m(a^\dag)^{n-2} \\
&=\cdots=(a^\dag+m\lambda)^n a^m.
\end{split}
\end{equation}

%10
Therefore, by \eqref{eq20}, \eqref{eq21} and \eqref{eq22}, we obtain the following theorem.

\begin{theorem}

For $m,n \in \mathbb{N}$ and $\lambda\neq0$, we have in $S_\lambda$ the normal orderings given by
\begin{equation*}
\begin{split}
&a^m a^\dag=(a^\dag+m \lambda)a^m,\quad a(a^\dag)^n=(a^\dag+\lambda)^na, \\
&a^m(a^\dag)^n=(a^\dag+m\lambda)^n a^m.
\end{split}
\end{equation*}

\end{theorem}

\medskip

Now, we observe from Theorem 3 that
\begin{equation}\label{eq23}
\begin{split}
(a^\dag a)^{n+1}&=(a^\dag a)(a^\dag a)^n=a^\dag a\sum_{k=0}^n {n \brack k}_\lambda (a^\dag)^ka^n \\
&=a^\dag \sum_{k=0}^n {n \brack k}_\lambda a(a^\dag)^ka^n=a^\dag\sum_{k=0}^n {n \brack k}_\lambda (\lambda+a^\dag)^k aa^n \\\
&=a^\dag \sum_{k=0}^n {n \brack k}_\lambda \sum_{j=0}^k \binom{k}{j}\lambda^{k-j}(a^\dag)^ja^{n+1} \\
&=\sum_{j=0}^n\sum_{k=j}^n {n \brack k}_\lambda \binom{k}{j}\lambda^{k-j} (a^\dag)^{j+1}a^{n+1} \\
&=\sum_{j=1}^{n+1} \bigg(\sum_{k=j-1}^n {n \brack k}_\lambda \binom{k}{j-1}\lambda^{k+1-j}\bigg)(a^\dag)^ja^{n+1}.
\end{split}
\end{equation}

%11
On the other hand, by Theorem 1, we get
\begin{equation}\label{eq24}
\begin{split}
(a^\dag a)^{n+1}=\sum_{j=0}^{n+1}{n+1 \brack j}_\lambda(a^\dag)^j a^{n+1}=\sum_{j=1}^{n+1}{n+1 \brack j}_\lambda(a^\dag)^ja^{n+1}.
\end{split}
\end{equation}

Therefore, by \eqref{eq23} and \eqref{eq24}, we obtain the following theorem.
\begin{theorem}

Let $n,j\in\mathbb{Z}$ with $n\geq0$ and $j\geq1$. In $S_\lambda$, the unsigned $\lambda$-Stirling numbers of the first kind satisfy the following recurrence relation:
\begin{equation*}
\begin{split}
{n+1 \brack j}_\lambda=\sum_{k=j-1}^n{n \brack k}_\lambda\binom{k}{j-1}\lambda^{k+1-j}={n \brack j-1}_\lambda+\sum_{k=j}^n{n \brack k}_\lambda\binom{k}{j-1}\lambda^{k+1-j}.
\end{split}
\end{equation*}

\end{theorem}

\medskip

For $n\geq1$, by \eqref{eq16} and \eqref{eq18}, we have the $\lambda$-analogues of Boole's relations in the $\lambda$-Shift algebra given by
%12
\begin{equation*}
\begin{split}
(a^\dag a)^n=\langle a^\dag \rangle_{n,\lambda} a^n,\quad ((a^\dag+r)a)^n=\langle a^\dag+r\rangle_{n,\lambda}a^n.
\end{split}
\end{equation*}

Now, we define the $\lambda$-analogues of $n!$ as (see \eqref{eq08})
\begin{equation}\label{eq25}
\begin{split}
(0)_\lambda !=1,\quad (n)_\lambda !=\langle 1\rangle_{n,\lambda}=\sum_{k=0}^n{n \brack k}_\lambda,\quad (n\geq1).
\end{split}
\end{equation}
Note that $\displaystyle \lim_{\lambda\rightarrow 1}(n)_\lambda!=n!$.

From \eqref{eq08} and \eqref{eq09}, we note that
\begin{equation}\label{eq26}
\begin{split}
\frac{1}{k!}\bigg(-\frac{\log(1-\lambda t)}{\lambda}\bigg)^k=\sum_{n=k}^\infty{n \brack k}_\lambda\frac{t^n}{n!},
\end{split}
\end{equation}

and
\begin{equation}\label{eq27}
\begin{split}
\frac{1}{k!}\bigg(-\frac{\log(1-\lambda t)}{\lambda}\bigg)^k\bigg(\frac{1}{1-\lambda t}\bigg)^{\frac{r}{\lambda}}=\sum_{n=k}^\infty{n+r \brack k+r}_{r,\lambda}\frac{t^n}{n!},
\end{split}
\end{equation}
where $k$ is nonnegative integer.

Thus, by \eqref{eq26} and \eqref{eq27}, we get
\begin{equation}\label{eq28}
\begin{split}
\sum_{n=k}^\infty{n+1 \brack k+1}_{1,\lambda}\frac{t^n}{n!}&=\bigg(\frac{1}{1-\lambda t}\bigg)^{\frac{1}{\lambda}}\frac{1}{k!}\bigg(-\frac{1}{\lambda}\log(1-\lambda t)\bigg)^k\\
&=\sum_{l=0}^\infty \frac{\langle 1\rangle_{l,\lambda}}{l!}\ t^l\sum_{m=k}^\infty{m \brack k}_\lambda \frac{t^m}{m!}\\
&=\sum_{n=k}^\infty \bigg(\sum_{l=0}^{n-k}\binom{n}{l}\langle 1 \rangle_{l,\lambda}{n-l \brack k}_\lambda\bigg)\frac{t^n}{n!}.
\end{split}
\end{equation}

%13
Comparing the coefficients on both sides of \eqref{eq28}, we get
\begin{equation}\label{eq29}
\begin{split}
{n+1 \brack k+1}_{1,\lambda}=\sum_{l=0}^{n-k}\binom{n}{l}\langle 1\rangle_{l,\lambda}{n-l \brack k}_\lambda.
\end{split}
\end{equation}

From \eqref{eq09} and \eqref{eq29}, we note that
\begin{equation}\label{eq30}
\begin{split}
\langle 1+m\lambda\rangle_{n,\lambda}&=\sum_{k=0}^n{n+1 \brack k+1}_{1,\lambda}m^k\lambda^k=\sum_{k=0}^n\sum_{l=0}^{n-k}\binom{n}{l}\langle 1\rangle_{l,\lambda}{n-l \brack k}_\lambda m^k\lambda^k\\
&=\sum_{l=0}^n\bigg(\sum_{k=0}^{n-l}{n-l \brack k}_\lambda m^k \lambda^k\bigg)\binom{n}{l}\langle 1 \rangle_{l,\lambda}=\sum_{l=0}^n\binom{n}{l}\langle 1\rangle_{l,\lambda}\langle m\lambda\rangle_{n-l,\lambda}.
\end{split}
\end{equation}

Therefore, by \eqref{eq30}, we obtain the following theorem.

\begin{theorem}

For $n\geq0$, we have
\begin{equation*}
\begin{split}
\langle 1+m\lambda\rangle_{n,\lambda}=\sum_{l=0}^n\binom{n}{l}\langle 1\rangle_{l,\lambda}\langle m\lambda\rangle_{n-l,\lambda}.
\end{split}
\end{equation*}

\end{theorem}

\medskip

%14
We note that
\begin{equation*}
\begin{split}
\langle 1+m \rangle_n=\lim_{\lambda\rightarrow1}\langle 1+m\lambda\rangle_{n,\lambda}=\sum_{l=0}^n\binom{n}{l}l!\langle m\rangle_{n-l}.
\end{split}
\end{equation*}

For any $\alpha\in\mathbb{R}$, we define the $\lambda$-shift operator $\delta_\lambda^{\alpha}$ by
\begin{equation}\label{eq31}
\begin{split}
\delta_\lambda^\alpha f(x)=f(x+\lambda\alpha).
\end{split}
\end{equation}
Then we see that
\begin{equation*}
\begin{split}
\delta_\lambda x-x\delta_\lambda=\lambda\delta_\lambda,
\end{split}
\end{equation*}
where $\delta_\lambda=\delta_\lambda^1$, and $x$ denotes the `multiplication by $x$' operator.

In the $\lambda$-shift algebra $S_{\lambda}$, a concrete representation is given by the operators $a^\dag\mapsto x$ and $a\mapsto \delta_\lambda$.
From Theorems 1 and 3, we note that
\begin{equation*}
\begin{split}
\delta_\lambda^r x^s=(x+r\lambda)^s\delta_\lambda^r,
\end{split}
\end{equation*}

and
\begin{equation}\label{eq33}
\begin{split}
(x\delta_\lambda)^n=\sum_{k=0}^n{n \brack k}_\lambda x^k\delta_\lambda^n,\quad (n\geq0).
\end{split}
\end{equation}

Now, we observe from \eqref{eq31} and \eqref{eq33} that
\begin{equation}\label{eq34}
\begin{split}
(x\delta_\lambda)^ne^x=\sum_{k=0}^n{n \brack k}_\lambda x^k\delta_\lambda^ne^x=\sum_{k=0}^n{n \brack k}_\lambda x^k e^{(n\lambda+x)}.
\end{split}
\end{equation}

%15
By \eqref{eq08} and \eqref{eq34}, we get
\begin{equation}\label{eq35}
\begin{split}
e^{-(n\lambda+x)}(x\delta_\lambda)^ne^x=\langle x \rangle_{n,\lambda}=\sum_{k=0}^n{n \brack k}_\lambda x^k.
\end{split}
\end{equation}

In particular, for $x=1$, we have (see \eqref{eq25})
\begin{equation}\label{eq36}
\begin{split}
e^{-(n\lambda+x)}(x\delta_\lambda)^ne^x \Big\rvert_{x=1}=\sum_{k=0}^n{n \brack k}_\lambda=(n)_\lambda!.
\end{split}
\end{equation}

Therefore, by \eqref{eq35} and \eqref{eq36}, we obtain the following theorem.
\begin{theorem}

For $n\in\mathbb{N}$, in $S_\lambda$, we have
\begin{equation*}
\begin{split}
e^{-(x+n\lambda)}(x\delta_\lambda)^ne^x=\langle x \rangle_{n, \lambda}=\sum_{k=0}^n{n \brack k}_\lambda x^k.
\end{split}
\end{equation*}

In particular, for $x=1$, we get
\begin{equation*}
\begin{split}
e^{-(x+n\lambda)}(x\delta_\lambda)^ne^x\Big\rvert_{x=1}=\sum_{k=0}^n{n \brack k}_\lambda=(n)_\lambda!.
\end{split}
\end{equation*}
\end{theorem}

\medskip

From Theorem 6, we note that
%16
\begin{equation}\label{eq37}
\begin{split}
e^{-(x+(m+n)\lambda)}(x\delta_\lambda)^{m+n}e^x\Big\rvert_{x=1}=(m+n)_\lambda!.
\end{split}
\end{equation}

On the other hand, by \eqref{eq33}, we get
\begin{equation}\label{eq38}
\begin{split}
(x\delta_\lambda)^{m+n}&=(x\delta_\lambda)^m(x\delta_\lambda)^n=\sum_{j=0}^m{m \brack j}_\lambda \sum_{k=0}^n{n \brack k}_\lambda x^j\delta_\lambda^m x^k\delta_\lambda^n\\
&=\sum_{j=0}^m{m \brack j}_\lambda \sum_{k=0}^n{n \brack k}_\lambda x^j(x+m\lambda)^k\delta_\lambda^{m+n}.
\end{split}
\end{equation}

From \eqref{eq38}, we have
\begin{equation}\label{eq39}
\begin{split}
(x\delta_\lambda)^{m+n}e^x&=\sum_{j=0}^m{m \brack j}_\lambda \sum_{k=0}^n{n \brack k}_\lambda x^j(x+m\lambda)^k\delta_\lambda^{m+n}e^x\\
&= \sum_{j=0}^m{m \brack j}_\lambda \sum_{k=0}^n{n \brack k}_\lambda x^j(x+m\lambda)^k e^{x+(m+n)\lambda}.
\end{split}
\end{equation}

Thus, by \eqref{eq39}, we get
\begin{equation}\label{eq40}
\begin{split}
e^{-(x+(m+n)\lambda)}(x\delta_\lambda)^{m+n}e^x=\sum_{j=0}^m{m \brack j}_\lambda \sum_{k=0}^n{n \brack k}_\lambda x^j (x+m\lambda)^k.
\end{split}
\end{equation}

From \eqref{eq37} and \eqref{eq40}, we have

\begin{equation}\label{eq41}
\begin{split}
(m+n)_\lambda!&=e^{-(x+(m+n)\lambda)}(x\delta_\lambda)^{n+m}e^{x}\Big\rvert_{x=1}\\
&=\sum_{j=0}^m{m \brack j}_\lambda \sum_{k=0}^n{n \brack k}_\lambda (1+m\lambda)^k\\
&=\sum_{j=0}^m{m \brack j}_\lambda\langle 1+m\lambda\rangle_{n,\lambda}.
\end{split}
\end{equation}

%17
By Theorem 5 and \eqref{eq41}, we get
\begin{equation}\label{eq42}
\begin{split}
(m+n)_\lambda!&=\sum_{j=0}^m{m \brack j}_\lambda\langle 1+m\lambda \rangle_{n,\lambda}\\
&=\sum_{j=0}^m{m \brack j}_\lambda\sum_{k=0}^n \binom{n}{k} \langle m\lambda\rangle_{n-k,\lambda}\langle 1\rangle_{k,\lambda}.
\end{split}
\end{equation}

Therefore, by \eqref{eq42}, we obtain the following theorem.
\begin{theorem}
For $m,n\geq0$, we have
\begin{equation*}
\begin{split}
(m+n)_\lambda!=\sum_{j=0}^m{m \brack j}_\lambda\langle 1+m\lambda \rangle_{n,\lambda}=\sum_{j=0}^m\sum_{k=0}^n{m \brack j}_\lambda \binom{n}{k} \langle m\lambda \rangle_{n-k,\lambda}\langle 1\rangle_{k,\lambda}.
\end{split}
\end{equation*}
\end{theorem}

\begin{remark}
(a) Taking the limit as $\lambda \rightarrow 1$, we see from Theorem 7 that
\begin{equation*}
\begin{split}
(m+n)!=\lim_{\lambda\rightarrow1}(m+n)_\lambda!=\sum_{j=0}^m\sum_{k=0}^n{m \brack j} \binom{n}{k} \langle m \rangle_{n-k}k!,\quad (\text{see \cite{16}}).
\end{split}
\end{equation*}
This was discovered by Mez\H{o} in \cite{16} which is dual to Spivey's identity (see \cite{20}, \cite{21}). \par
(b) As the identities in Theorem 7 are obviously symmetric in $m$ and $n$, we get the following symmetric identities:
\begin{equation*}
\begin{split}
&\sum_{j=0}^m{m \brack j}_\lambda\langle 1+m\lambda \rangle_{n,\lambda}=\sum_{j=0}^n{n \brack j}_\lambda\langle 1+n\lambda \rangle_{m,\lambda}, \\
&\sum_{j=0}^m\sum_{k=0}^n{m \brack j}_\lambda \binom{n}{k} \langle m\lambda \rangle_{n-k,\lambda}\langle 1\rangle_{k,\lambda}
=\sum_{j=0}^n\sum_{k=0}^m{n \brack j}_\lambda \binom{m}{k} \langle n\lambda \rangle_{m-k,\lambda}\langle 1\rangle_{k,\lambda}.
\end{split}
\end{equation*}
\end{remark}

\section{Conclusion}
In this paper, we introduced the $\lambda$-shift algebra $S_{\lambda}$ as a $\lambda$-analogue of the shift algebra $S$ and derived normal ordering results in $S_{\lambda}$ where the unsigned $\lambda$-Stirling numbers of the first kind and
the $\lambda$-$r$-Stirling numbers of the first kind appear as the coefficients. In addition, from these normal ordering results we derived some properties about the unsigned $\lambda$-Stirling numbers of the first kind. \par
There are various methods that can be used in order to find some results on special numbers and polynomials. These include generating functions, combinatorial methods, umbral calculus, $p$-adic analysis, differential equations, analytic number theory, probability and statistics, operator theory, special functions and mathematical physics. \par
It is one of our future projects to continue to explore many special numbers and polynomials with these various tools.

%\begin{equation}\label{eq01}
%\begin{split}
%
%\end{split}
%\end{equation}

%\begin{equation}\label{eq01}
%\begin{split}
%
%\end{split}
%\end{equation}

%\begin{equation}\label{eq01}
%\begin{split}
%
%\end{split}
%\end{equation}

%\begin{equation}\label{eq01}
%\begin{split}
%
%\end{split}
%\end{equation}

%\begin{equation}\label{eq01}
%\begin{split}
%
%\end{split}
%\end{equation}

%\begin{equation}\label{eq01}
%\begin{split}
%
%\end{split}
%\end{equation}

%\begin{equation}\label{eq01}
%\begin{split}
%
%\end{split}
%\end{equation}

\bigskip

%\noindent{\bf{Acknowledgments}} \\
%The author would like to thank the referees for the detailed and valuable comments that helped
%improve the original manuscript in its present form.% Also,
%%
%The authors thank Jangjeon Institute for Mathematical Science for the support of this research.

\vspace{0.1in}

\noindent{\bf {Availability of data and material}} \\
Not applicable.

\vspace{0.1in}

%\noindent{\bf{Funding}} \\
%This work was supported by the Basic Science Research Program, the National
  %             Research Foundation of Korea,
  %             (NRF-2021R1F1A1050151).
\vspace{0.1in}

\noindent{\bf{Ethics approval and consent to participate}} \\
The authors declare that there is no ethical problem in the production of this paper.

\vspace{0.1in}

\noindent{\bf {Competing interests}} \\
The authors declare no conflict of interest.

\vspace{0.1in}

%\noindent{\bf{Consent for publication}} \\
%The authors want to publish this paper in this journal.

\vspace{0.1in}

%\noindent{\bf{Author' Contributions}}\\
%HKK  structured and wrote the whole paper. DSL and HKK checked the results of the paper and
%completed the revision of the article. All authors read and approved the final manuscript.
%\vspace{0.1in}
%
%\noindent{\bf{Author details}}\\
%Department of Mathematics Education, Daegu Catholic University, Gyeongsan 38430, Republic of Korea.

%\vspace{0.1in}
\

\end{document}